\documentclass[12pt,twoside,a4paper]{article}
\usepackage[pdftex]{hyperref}
\usepackage{amssymb}
\usepackage[centertags]{amsmath}

\usepackage{amsthm}
\usepackage{graphicx}
\font\teneufm=eufm10
\font\seveneufm=eufm7
\font\fiveeufm=eufm5
\newfam\eufmfam
\textfont\eufmfam=\teneufm
\scriptfont\eufmfam=\seveneufm
\scriptscriptfont\eufmfam=\fiveeufm

\usepackage{xcolor}

\newcommand{\blue}[1]{\textcolor{blue}{ {#1}}}
\newcommand{\violet}[1]{\textcolor{violet}{ {#1}}}

\newtheorem{thm}{Theorem}

\newtheorem{defn}{Definition}
\newtheorem{pro}{Proposition}
\newtheorem{lem}{Lemma}
\newtheorem{rem}{Remark}
\newtheorem*{rem*}{Remark}
\newtheorem{cor}{Corollary}

\renewcommand{\proof}{\textsc{proof. }}
\newcommand{\eproof}{\hfill $\Box$}

\renewcommand{\o}{\omega}

\newcommand{\e}{\varepsilon}
\renewcommand{\a}{\alpha}
\renewcommand{\b}{\beta}
\renewcommand{\d}{\delta}

\newcommand{\noi}{\noindent}

\renewcommand{\l}{\lambda}
\newcommand{\s}{\sigma}

\newcommand{\z}{\zeta}
\renewcommand{\t}{\tau}
\newcommand{\g}{\gamma}

\newcommand{\N}{\mathbb{N}}

\newcommand{\Z}{\mathbb{Z}}
\newcommand{\R}{\mathbb{R}}
\newcommand{\C}{\mathbb{C}}
\renewcommand{\th}{{\mathtt{h}}}
\newcommand{\tw}{{\mathtt{w}}}
\newcommand{\tr}{{\mathtt{r}}}
\newcommand{\set}[1]{{\left\{#1\right\}}}
\newcommand{\pa}[1]{{\left(#1\right)}}
\newcommand{\abs}[1]{{\left|#1\right|}}
\newcommand{\norm}[1]{{\left |#1\right |}}
\newcommand{\im}{{\rm i}}
\newcommand{\jap}[1]{| #1 |}
\newcommand{\und}[1]{\underline{#1}}
\newcommand{\cH}{{\mathcal H}}
\newcommand{\cccp}{{\mathfrak c}}

\makeindex

\title{A quantitative Birkhoff Normal Form \\for the hinged-hinged beam equation
\\
with geometric nonlinearity}

\author{L. Di Gregorio, W. Lacarbonara}

\begin{document}

\allowdisplaybreaks
\maketitle

{\abstract{
We consider an undamped nonlinear hinged-hinged beam 
with stretching nonlinearity as an infinite dimensional 
hamiltonian system.
 We obtain
analytically a quantitative  
Birkhoff Normal Form, via a nonlinear 
coordinate transformation that allows us to integrate
the system up to a small reminder, providing a very precise description
of small amplitude solutions over large time scales. 
The optimization of the  involved estimates  yields results obtained for realistic values of the physical quantities and of the perturbation parameter.
}}

{\small 
\tableofcontents
}

\noindent
\section{Introduction and main results}

We consider the  
dimensionless
nonlinear  beam equation
with stretching nonlinearity
\begin{equation}\label{checco}
\displaystyle
u_{tt}+u_{xxxx}- 
\left(m+ {\frac 1{2\pi}\int_0^{\pi}u_x^2\, dx}\right)u_{xx}
=0\,,
\end{equation}
for $t\in\mathbb R$ and   $x\in[0,\pi]$,
with the following hinged-hinged boundary conditions:
\begin{equation}\label{bc}
u(t,0)=u(t,\pi)=u_{xx}(t,0)=u_{xx}(t,\pi)=0
\end{equation}
and initial data
\begin{equation}\label{pisolino}
u_0(x):=u(0,x)\,,\qquad v_0(x):=u_t(0,x)\,.
\end{equation}
Here 
$\sqrt{I/A}u$ is the vertical displacement
and $m=\frac{L^2 P}{\pi^2 EI}$ indicates the nondimensional axial force,
where
 $L,I,A,E,P$ are, respectively,
  the  length of the beam,
the moment of inertia,
the cross-section area,
the Young modulus, and
the tensile axial force (possibly also negative entailing compressive force); see Subsection \ref{zabaione} below
for the derivation of \eqref{checco}
from the dimensional physical  equation.
\\
Being conservative, the equation in \eqref{checco}
has an infinite dimensional hamiltonian structure. Indeed, 
for $j\in\mathbb N:=\{1,2,3,\ldots\},$
by
letting
$\omega_j^2:=j^4+mj^2$ and $\phi_j(x)
:=\sqrt{2/ \pi}\sin jx$, respectively, denote the eigenvalues and the eigenfunctions
of the Sturm--Liouville operator 
$(\partial_{xxxx}-m\partial_{xx})$ 
on $[0,\pi]$,
the Hamiltonian can be expressed as
\begin{equation}\label{tommy}
H(p,q)=
\sum_{j\in \N} \omega_j I_j+\frac 1{8\pi} \left(\sum_{j\in \N} \frac{j^2}{\omega_j}q_j^2\right)^2\,,\qquad
\omega_j:=\sqrt{j^4+mj^2}\,,\qquad
I_j:=\frac{p_j^2+q_j^2}{2}\,,
\end{equation}
with
$q=(q_1,q_2,\ldots,)$, $p=(p_1,p_2,\ldots)$ spanning 
a suitable Hilbert space of sequencies
(typically a $L^2$-based Sobolev space).  
Then, given a sufficiently smooth solution $t\to \big(p(t), q(t)\big)$ of
Hamilton's equations of $H$,
one finds that 
 $u(t,x):=
\sum_{j\geq 1}\frac{q_j(t)}{\sqrt{\omega_j}}\phi_j(x)$
is a  
solution of \eqref{checco}.
Assuming $m>-1$ we have that $\omega_j>0$
and
the origin is an elliptic equilibrium for 
$H$. Then, a classical way to describe the behaviour
of small amplitude solutions and to obtain 
information on their (long) time stability is
by the Birkhoff Normal Form (BNF).
Indeed if the linear frequencies 
$\omega_j$ are  non-resonant
up to order $2d>0$, namely 
$\omega\cdot k$
is sufficiently bounded away from zero
for all $k\in\mathbb Z^{\mathbb N}$
with $|k|\leq 2d$, after a
close-to-the-identity canonical 
change of variables, in a sufficiently 
small neighbourhood    of  the origin
the Hamiltonian is in
BNF  up to  order $2d$, namely 
$H=\mathcal N(I)+R( p,q)
$, where $\mathcal N(I)$
is a suitable  polynomial  of degree $d$ in the
``actions''
$I=(I_j)_{j\in\mathbb N}$  and
 $R=O(|(p,q)|^{2d+2} )$.  
Since the nonlinear term $\mathcal N$
is \textsl{integrable}, 
the BNF allows a precise 
description of the  solutions 
with initial data $\big(p(0),q(0)\big)$ satisfying 
$\epsilon:=2^{-1/2}\sqrt{|p(0)|^2+|q(0)|^2}\leq \epsilon_0$
   up to  times  
$|t|\leq {\rm T}_0\epsilon^{-2d}$, for suitably small
$\epsilon_0$ and ${\rm T}_0$.
This immediately reads as
 a time stability result for 
  \eqref{checco}
with $\epsilon$-small initial data
$u_0(x)$ and $v_0(x)$.

\bigskip

The problem of long time stability of PDEs
close to an elliptic equilibrium via BNF
 has been widely studied 
in particular for finite-regularity (Sobolev)
 initial data.
 The first polynomial bounds on the stability time
 for hamiltonian 
 PDEs was showed in 
the seminal works \cite{Bambusi:2003} and \cite{Bambusi-Grebert:2006}.
   More recently, similar results have been obtained also for PDEs with nonlinearities containing derivatives, see \cite{Delort-2009},
   \cite{Yuan-Zhang},\cite{Delort-2015},\cite{Berti-Delort},\cite{FI}. 
All the above  results
rely on suitable non resonance conditions
on $\omega$ avoiding exact internal resonances
and 
ensuring suitable lower bounds on the quantity
 $|\omega\cdot k|$.
  These arithmetic conditions are typically achieved by suitably tuning some 
 ``parameter'' which modulates the linear frequencies $\omega_j$;
 the parameter can be either ``internal", such as 
 $m$ in \eqref{checco} or
 the capillarity of the fluid in the case of water waves, or ``external", as, e.g., convolutions and multiplicative potentials for the nonlinear Schr\"odinger equation. 
 The non resonance condition is imposed 
 for ``most'' (in a measure sense) values
 of the parameter(s).
 We note that in two or higher space dimension
 could be not possible to impose such conditions
 at any order,  see, e.g., \cite{I}, \cite{BFGI}.
\\
 All the above results  regard polynomial stability times in Sobolev spaces. 
Passing from polynomial estimates to exponential-type ones 
is not obvious and can be achieved
requiring  more regularity  on the initial data,
see, e.g.,   
 \cite{Faou-Grebert:2013}, \cite{Cong},\cite{BMP} and \cite{FM}.

\smallskip

Recollecting we can say that, 
from the mathematical point of view,
all the efforts in performing BNF for PDEs
by methods of Perturbation Theory\footnote{Mainly
in the hamiltonian or reversible case.}
have been in the direction of extending time stability
and/or the range of applicability to larger and larger classes
of equations (for example containing derivatives in the nonlinearitiy).
Almost no attention was given to the problem 
of how small must be the perturbative 
parameter, which measures the amplitude
of initial data. 
Indeed the ``applicability threshold'',
giving an upper bound on the maximal size of 
 initial data, is, in general, not even evaluated\footnote{With the remarkable
exception of \cite{BMP} where, however, it is extremely
small .....}
but, surely, extremely
small excluding any possibility of an application 
to concrete cases from Physics
and/or Engineering.
Actually this is a well known problem in Perturbation 
Theory already in finite dimension where, instead,
some results exist (for example in Celestial Mechanics).
However, as we said above, almost nothing 
is known in the infinite dimensional case of PDEs.

On the contrary  in this paper we put  all our attention
in evaluating the applicability threshold, 
optimising  the methods in order to obtain
the larger possible value. Analogously 
we also evaluate the stability time not simply its 
``order'' as a power of the perturbative parameter.
See, in particular, formulas 
\eqref{trescul0}, \eqref{burrata0}, 
\eqref{e1} and \eqref{burrataBIS} below.
It turns out that both the applicability threshold
and the stability time are realistic and fit the physical
parameters; see figure 1 below.
As far as we know, \emph{this is the first 
 result
of this kind in hamiltonian perturbation theory for 
PDEs}.

\subsection{Main results}

We look for solutions $u(t,x)$ of equations 
\eqref{checco}--\eqref{ginger}
such that, for every fixed time $t$, the function
$x\to u(t,x)$ belongs to 
the Hilbert space
 $H^s_*\subset H^s$, $s\geq 0$, 
 of functions 
 $u(x)=\sum_{j\in\N}u_j\phi_j(x)$, 
 $u_j\in\mathbb R$,
 with finite norm
 $$
|u|_s:=
\Big(\sum_{j\in\N}u_j^2 j^{2s}\Big)^{1/2}
=\Big(\int_0^\pi |\partial_x^s u(x)|^2dx
\Big)^{1/2}\,,$$
where
$\partial_x^s u(x)$
is the (Sobolev) 
weak
$s$--derivative of $u(x)$.
Note that
  $\{\phi_j(x)\}_{j\in\mathbb N}$ 
  is an orthonormal basis of the space $H^s_*$
 endowed with the scalar product
 $\langle u, v\rangle:=
 \sum_{j\in\N}j^{2s}u_j v_j$.
 Note that, for $s>1/2$, 
all the functions $u\in H^s_*$ are continuous with
\begin{equation}\label{mac} 
\max_{x\in[0,\pi]}|u(x)|\leq c_s |u|_s\,,\qquad
c_s^2:=\frac{2}{\pi}\sum_{j\in\N}j^{-2s}
\end{equation}
and satisfy $u(0)=u(\pi)=0$.
Moreover, for $s>5/2$, $u\in H^s_*\subset \mathcal C^2$
satisfies  $u(0)=u(\pi)=u_{xx}(0)=u_{xx}(\pi)=0$.
\\
Given 
 an open interval of times $\mathcal I\ni 0$,
 $s>1/2$ and two initial data
 $u_0\in H^{s+2}_*$, $v_0\in H^s_*$,
we say that a $C^2$ function 
$t\mapsto u(t,\cdot)\in H^{s+2}_*$
is a {\sl weak solution} of 
\eqref{checco}-\eqref{ginger}
if it satisfies \eqref{ginger}
and for every test function $\eta\in\mathcal C^\infty_c((0,\pi))$
we have
$$
\int_0^{\pi}\left[u_{tt}\eta+u_{xx}\eta_{xx}+\left(m+\frac 1{2\pi}\int_0^{\pi}u_x^2\right) u_x \eta_x\right]dx=0\,,
\qquad
\forall\, t\in \mathcal I\,.
$$
We now state our two main stability results:
the former follows by the fourth order BNF,
the latter  by the sixth order one.

\begin{thm}\label{patata}
Let $m\geq -1/2$, $s>1/2$ and assume that
$u_0\in H^{s+2}_*$ and 
$v_0\in H^{s}_*$
satisfy 
 \begin{equation}\label{trescul0}
\e:=\max\{|u_0|_1\,,\ 
|v_0|_0/(1+m)\}\ \leq \e_0:=
0.08\sqrt\mu\,,\qquad
\mu:=
\textstyle
\min\left\{2\sqrt{1+m},\frac32 \sqrt{4+m}\right\}
 \,.
\end{equation}
Then the weak solution $u(t,x)$
of equation \eqref{checco}-\eqref{ginger}
satisfies 
\begin{equation}
 \sup_{x\in[0,\pi]}|u(t,x)| 
  \leq 2\e\,,\ \ 
\forall\, |t|\leq T_0\mu\e^{-4}\,,\ \ 
T_0\geq
0.008(1-68 \e^2/\mu)\geq
 0.0045
\,.
\label{burrata0}
\end{equation}
\end{thm}

\begin{rem}\label{alice}
An important fact both from an abstract
point of view  and for concrete applications
is that  on the initial data 
$u_0$, resp. $v_0$,
we require smallness conditions
(see trescul)
\emph{only in the very weak norms} 
$|\cdot|_1$, resp. $|\cdot|_0$.
 A  slightly stronger formulation of the above
result is given in Theorem
\ref{patata2} below where, in particular, 
we only require that the very weak norm 
$|\cdot|_{-1}$ of $v_0$ is small (see formula \eqref{trescul03} below).
\end{rem}

Quite better estimates hold making suitable restriction
on $m$ allowing no sixth order resonances.
First set
\begin{equation}\label{e1}
\e_1:=\min\left\{\e_0, \sqrt[4]{\frac{|m| \mu}{ 3222\cccp^2}}
\right\}\,,\qquad
T_1:=
\frac{|m|\mu}{6595\cccp^3
\Big(
1+30\frac{\cccp}{|m|\mu}\e^2+
90\frac{\cccp}{\mu}\e^2
\Big)}\,,
\end{equation}
where
\begin{equation}\label{cccp}
\cccp
:=\left\{
\begin{array}{ll}
\displaystyle 1\ \ \ \ &\mathrm{if}\ \ \ m\geq 0\,,
\\
\displaystyle (1+m)^{-1/2}\ \ \ \ &\mathrm{if}\ \ \ 
-1/2\leq m<0\,.
\end{array}
\right.
\end{equation}
For example for $m=1$ we have
$\e_1=\e_0\sim 0.13$, $T_1\geq 1.9 \times 10^{-4}$.

\begin{thm}\label{patataBIS}
In addition to the assumptions of Theorem
\ref{patata}
 we require that
 $-1/2\leq m\leq 1$, $m\neq 0$ and $\e\leq \e_1$.
Then
\begin{equation}
 \sup_{x\in[0,\pi]}|u(t,x)| 
  \leq 2.1\e\,,\qquad 
\forall\, |t|\leq T_1\e^{-6}\,.
\label{burrataBIS0}
\end{equation}
\end{thm}
\noi
A  slightly stronger version of the above result is given in Theorem \ref{patataTER}
below.

\bigskip

\textbf{Idea of the proof}
\\
In the landscape of 
(hamiltonian)
Perturbation Theory of PDEs 
 the very effective estimates \eqref{trescul0}-\eqref{burrataBIS} might seem quite surprising.
 The reason they hold true is twofold: 
 (i) the very 
 peculiar structure of the nonlinearity in \eqref{checco}
 and
 (ii)
  the technical machinery we use that allows us 
 to work in very low (spatial) regularity.
 
 (i)
 The
 ``Kirchhoff type'' nonlinearity 
 $-\frac1{2\pi}u_{xx}\int_0^{\pi}u_x^2\, dx$
 generates,
 in the Hamiltonian in \eqref{tommy}, 
 the  quartic term
 $\frac 1{8\pi} (\sum_{j\in \N} \frac{j^2}{\omega_j}q_j^2)^2$ .
 Passing to complex variables through
 the symplectic transformation 
 $z_j=\frac1{\sqrt 2}(q_j+\mathrm ip_j)$
 (see \eqref{tommaso} below),
 the quadratic part becomes
  $\frac{1}{32\pi}\big( \sum_{j\in \N} \frac{j^2}{\omega_j}(z_j+\bar z_j)^2\big)^2$, which has the special property
  that every monomial contains the index $k$ in the form
  $z_k^{\alpha_k} \bar z_k^{\beta_k}$
  with $\alpha_k+\beta_k$ \textsl{even}
  (see \eqref{crostata}).
  This fact greatly simplifies the study
  of the non-resonance condition necessary
  for the fourth order BNF, which, in general,
corresponds to show that
$\omega_i\pm\omega_j\pm \omega_k\pm
\omega_\ell$ is  suitably bounded away from 
zero\footnote{Except for the cases 
$\omega_i-\omega_i+ \omega_k-
\omega_k$ et similia, which correspond to monomials
of the form $|z_i|^2 |z_k|^2=I_i I_k$,
that cannot be removed but are integrable terms.}.
 However, in the present case, due to the 
 above discussed  
 property of the nonlinearity, the non-resonance condition
 reduces to  consider only $2(\omega_i-\omega_j)$,
 which we prove being                                                                     uniformly bounded away from zero 
 for every $i>j\geq 1$
 (see Proposition \ref{pera} below).
 This means that \emph{there are no fourth
 order resonances at all}.
Then we can put the system in BNF up to order 
$2d=4$ (see Theorem \ref{sob})
and Theorem \ref{patata} follows.
\\
Note that the fourth order BNF procedure
preserves the special property above
(see formula \eqref{cacao}).
Then we can try to put the system 
in sixth  order BNF.
To do so, in general, one should consider
possible resonances of the type
$\omega_i\pm\omega_j\pm \omega_k\pm
\omega_\ell\pm\omega_m\pm \omega_n$,
which, in our case, reduce to $2(\omega_i\pm\omega_j\pm \omega_k)$
.
Obviously the quantity 
$\omega_i-\omega_j- \omega_k$
vanishes when
$m=0$ and $i,j,k$ is a Pythagorean triple, namely
$i^2=j^2+k^2$.
Then taking, let us say, 
$m\neq 0$ and
 $-1/2\leq m\leq 1$, we are able to prove  that
 $|\omega_i-\omega_j- \omega_k|$
 is uniformly bounded away from zero 
 for every $i>j>k\geq 1$, namely
\emph{there are even no sixth
 order resonances}
(see Proposition \ref{pera2}).
Then we can put the system in BNF up to order 
$2d=6$ (see Theorem \ref{sob6})
and Theorem \ref{patataBIS} follows.

(ii)
Recalling Remark \ref{alice}
we assume that the initial data in \eqref{pisolino}
satisfy
$u_0\in H^{s+2}_*$ and 
$v_0\in H^{s}_*$ for some
$s>1/2$.
However we only assume that the initial data
$u_0$ and $v_0$ 
are small in \emph{weaker norms},  
namely in the $H^1$ and $H^{-1}$ norms,
repectively
(see \eqref{trescul0}
or, better, \eqref{trescul03}).
This  greatly increases the applicability of our results
to real cases.
From a technical point of view, 
taking advantage of the 
machinery introduced in \cite{BMP}
(reminiscent of the Cauchy's majorants method),
for every $N>1$
we introduce an \emph{ad hoc} norm
$|\cdot|_{H^\sigma_N}$
 which is topologically equivalent to the 
 usual 
 one $|\cdot|_{H^\sigma}$,
 but on the first $N$ Fourier 
modes coincides with the norm $|\cdot|_{H^1}$,
see \eqref{enorme}. 
It is obvious that
for every fixed function $u$
we have that 
$|u|_{H^\sigma_N}$ tends to
$|u|_{H^1}$
 as $N\to\infty$.
 Then for every given initial 
 position\footnote{An analogous procedure 
 holds for $v_0$.}
 $u_0$, we take $N$ large enough
 so that $|u_0|_{H^{s+2}_N}$ is almost equal
 to $|u_0|_{H^1}$; so that
 imposing a smallness condition on
 $|u_0|_{H^1}$ is equivalent to impose 
 a smallness condition on
 $|u_0|_{H^{s+2}_N}$.
 Then we perform all the BNF procedure
 in the space $H^{s+2}$ endowed with the norm
 $|\cdot|_{H^{s+2}_N}$.
 The drawback is that we only that the solution
 $u(t,x)$ remains small only in the weaker norm
$|\cdot|_{H^{s+2}_N}$ (see \eqref{burrata}), while the standard norm
$|\cdot|_{H^{s+2}}$ could be very large.
However note that this implies a  bound 
on the $H^1$-norm and, therefore, on the sup-norm,
 namely the vertical displacement of the beam
 (see \eqref{burrata0}).

\begin{rem}\label{finitodim}
A very peculiar property of equation  \eqref{checco},
that we (apparently) do not use here, 
is the following:
 if the initial data in \eqref{pisolino}
  are trigonometric polynomials,
 then the solution of 
 is a trigonometric polynomial (in $x$) too.
 Actually the same holds true if the initial data 
 have Fourier support on whatever finite or infinite
 $\mathcal S\subseteq \mathbb N$.
Compare Remark \ref{pollastri}.
\end{rem}

\bigskip
The paper is divided as follows.
In the next subsection we discuss 
the application to the physical model.
In Section \ref{sec2} we describe the infinite dimensional 
hamiltonian structure of the beam equation \eqref{checco}.
In Section \ref{sec3} we study the fourth and sixth order
resonances. 
In Section \ref{sec4} we construct, using the estimates
on the resonances given in Section \ref{sec3}, 
the BNF up to order 4 and 6.
In Section \ref{sec5} we derive stability results
by the BNF.
 We prove some technical lemmata are proven in the Appendix.

\subsection{Application to the physical model}\label{zabaione}

Following  \cite{WK50} and \cite{D70} and \cite{MN95},
we consider a 
hinged beam of length $L$,
density  $\varrho$, cross section area $A$, 
   cross section 
 moment of inertia $I$ and
 Young's modulus $E$, under a stretching axial force $P$ due to the elongation.
It turns out that the deflection $w(\t,\xi)$ in a point $\xi\in[0,L]$
at time $\t\in\R$ is the solution
of the following PDE:
\begin{equation}\label{coblenza}
\varrho A w_{\t\t}+EI w_{\xi\xi\xi\xi}-\left(
P+\frac{EA}{2L}\int_0^L w_\xi^2 d\xi\right)
w_{\xi\xi}=0
\end{equation}
with hinged boundary conditions
\begin{equation}\label{betta2}
w(\t,0)=w(\t,L)=w_{\xi\xi}(\t,0)=w_{\xi\xi}(\t,L)=0
\end{equation}
and initial data
\begin{equation}\label{ginger2}
w(0,\xi)=w_0(\xi)\qquad
w_\t(0,\xi)=\mathtt v_0(\xi)\,.
\end{equation}


\noindent
Introducing dimensionless variables
\begin{equation}\label{cipolla}
t=\nu \t\,,\quad
\nu^2=\frac{\pi^4EI}{L^4\varrho A}\,,
 \quad x=\frac{\pi}{L}\xi\,,\quad 
m=\frac{L^2 P}{\pi^2 EI}\,,\quad
u(t,x)=\sqrt{\frac{A}{I}}w\left(\frac{t}{\nu},\frac{L}{\pi}x\right)\,,
\end{equation}
equation \eqref{coblenza} becomes
\eqref{checco}
with boundary conditions \eqref{bc}
and initial conditions
\begin{equation}\label{ginger}
u(0,x)=u_0(x)=:\sqrt{\frac{A}{I}}w_0\left(\frac{L}{\pi}x\right)\,,\qquad
u_t(0,x)=v_0(x)=:\frac{1}{\nu}\sqrt{\frac{A}{I}}\mathtt v_0\left(\frac{L}{\pi}x\right)\,.
\end{equation}

\begin{rem}
We choose the rescaling 
$u=\sqrt{A/I}w$ in \eqref{cipolla}
in order to have only one parameter, namely $m$,
in equation \eqref{checco}. A different (natural)
rescaling
would be $\tilde u=\frac{\pi}{L}w$, giving
 \begin{equation}\label{checcobis}
\tilde u_{tt}+\tilde u_{xxxx}- 
\left(m+ {\frac {\tilde m}{2}\int_0^{\pi}\tilde u_x^2\, dx}\right)\tilde u_{xx}
=0\,,\qquad
\tilde m:=\frac{AL^2}{\pi^3I}\,.
\end{equation}
Then $u$ and $\tilde u$ are related by
$$
\tilde u=\frac{\pi}{L}\sqrt{\frac{I}{A}}u\,.
$$
\end{rem}

\noi
As an example consider a beam with square 
 cross section of side $h$,
 so that
 the ratio between the
cross section moment of inertia 
and the cross section area is
$I/A=h^2/12$. 
We also take $h/L=10^{-2}$.
As initial data in \eqref{ginger2} we take $\mathtt v_0(\xi)=0$
and 
$$
w_0(\xi):=\d L p(\pi\xi/L)\,,\qquad \d>0
$$
$$
p(x):=-\frac{16}{5\pi^4}\left(x-\frac{\pi}{2}\right)^4
+
\frac{24}{5\pi^2}
\left(x-\frac{\pi}{2}\right)^2
-1\,.
$$
Note that $p$ is the polynomial of minimal degree which is even in $\pi/2$,
satisfies the hinged-hinged boundary conditions, namely
$p(0)=p(\pi)=p''(0)=p''(\pi)=0$,
and $\max_{[0,\pi]}|p|=1=-p(\pi/2)$;
moreover $|p|_{1}\sim 1.2583$. 
So the maximal initial displacement
is $\sup_{\xi\in[0,L]}|w_0(\xi)|=\d L$.
By \eqref{ginger} 
$u_0(x)=200\sqrt 3 \d p(x)$.
Then, recalling \eqref{trescul0},
$\e=|u_0|_1\sim 436\d$.
We choose $\d=10^{-4}$ and different values of $P$ such that 
$m=0,1$ or $-0.5$ so that
we have $T_0=0.0075$,  $T_1=7.8\times 10^{-4}$ or $T_1=2.9\times 10^{-5}$,
respectively.
Recalling \eqref{cipolla},
\eqref{burrata0} and \eqref{burrataBIS0},
the maximal displacement
$\sup_{\xi\in[0,L]}|w(\tau,\xi)|$
 remains
bounded by 
$2.1 \e/200\sqrt 3\sim 
2.643\times 10^{-4} L$
up to a time 
$|\tau|\leq T$
with
$T=2T_0\e^{-4}/\nu
$
in the case $m=0$
and
$T=T_1\e^{-6}/\nu
$
in the case $m=1$ and $-0.5$, respectively.
\bigskip

		 \begin{tabular}{|l|l|l|l|l|l|l|l|l|}
		 	\hline
				\text{\blue{Material}}& 
				$\blue{E}$ 
				$\scriptstyle{[10^{9}N \rm m^{-2}]}$
				& $\blue{\rho}$ $\scriptstyle{[kg\,\rm m^{-3}]}$& $\blue{L}\,\,\scriptstyle{[\rm m]}$ & $\blue{m}$ & \!$\blue{P}\,\, \scriptstyle{[kN]}$ & $				\blue{\nu}\,\scriptstyle{[s^{-1}]}$ & $\blue{T}\,\scriptstyle{[s]}$
				\\  
			 \hline 
			 	Steel & 200 & 7500  & 2 & 
				$0$ &\violet 0&  74 & \violet {56}  \\
		 	 	Steel & 200 & 7500  & 2 &$1$ &\violet{6.6}  & 74 &  \violet {1556}\\
				Steel & 200 & 7500  & 2 &$-0.5$ &\violet{-3.3}  & 74 &  \violet {58}\\
				Al 7075  & $70$   &2810  &  2 & $0$&\violet0&  71 & \violet {59}  \\
			 	Al 7075  & $70$   &2810  &  2 &$1$& \violet{2.3}&  71 &  \violet {1621} \\
				Al 7075  & $70$   &2810  &  2 &$-0.5$& \violet{-1.15}&  71 &  \violet {60} \\
 				 Al 7075  & $70$   &2810  
				 &  1 & $1$ &\violet{0.56}& 142 & \violet{811} \\
				Rubber &  0.004 & 1000  & 0.1 & $0$ &\violet0&  18 & \violet{232}  \\
			 	Rubber &  0.004 & 1000 & 0.1 & $1$ &\violet{$3 \!\cdot \!10^{-4}$} \!\!&  18 & \violet{6396} \\	
				Rubber &  0.004 & 1000 & 0.1 & $-0.5$ &\violet{$1.5 \!\cdot \!10^{-4}$} \!\!&  18 & \violet{238} \\	
			\hline
		\end{tabular}{\center{Fig. 1 Stability time for different values of the physical parameters}}

	\section{The infinite dimensional hamiltonian structure}
\label{sec2}
\noindent
Let us write \eqref{checco} as a 
 first order system
\begin{equation}
\left\{
\begin{array}{l}
\dot u=v\\
\dot v=-u_{xxxx}+\displaystyle
\left(m+ {\frac 1{2\pi}\int_0^{\pi}u_x^2\, dx}\right)u_{xx}
\end{array}
\right.
\end{equation}
with initial condition
\begin{equation}\label{gingerina}
u(0)=u_0(x)\in H^{s+2}_*\,,\qquad
v(0)=v_0(x)\in H^{s}_*
\end{equation}
in the phase space $H^{s+2}_* \times H^{s}_*([0,\pi])\ni (u,v)$. 
The associated Hamiltonian is
\begin{equation}
\mathtt H(u, v):=\int_{0}^{\pi}\left(\frac12 v^2+\frac 12 u_{xx}^2+\frac m2 u_{x}^2 \right)dx+
\frac{\pi}2\left(\displaystyle{\frac1{2\pi}\int_0^{\pi} u_{x}^2 dx}\right)^2\,.
\end{equation}
 Let us introduce Fourier coordinates $q=(q_1,q_2,\ldots,)$, $p=(p_1,p_2,\ldots)$, 
\begin{equation}\label{nduja}
u(x)=\sum_{j\in \N}\frac{q_j}{\sqrt{\omega_j}}\phi_j(x)\,\,\,\,\,\,\,\,\,\,\,\,\,\,v(x)=\sum_{j\in \N}p_j\sqrt{\omega_j}\phi_j(x)\,.
\end{equation}
The Hamiltonian reads (recall \eqref{tommy})
\begin{eqnarray*}
H(p,q)&=&\frac12\sum_{j\in \N}\omega_j p_j^2+\frac12\sum_{j\in \N}\frac{1}\omega_j j^4q_j^2+\frac m2 \sum_{j\in \N}\frac1\omega_j j^2 q_j^2 +\frac 1{8\pi} \left(\sum_{j\in \N} \frac{j^2}{\omega_j}q_j^2\right)^2\\
&=&
\sum_{j\in \N} \omega_j \frac{p_j^2+q_j^2}{2}+\frac 1{8\pi} \left(\sum_{j\in \N} \frac{j^2}{\omega_j}q_j^2\right)^2\,.
\end{eqnarray*}
The Hamilton equations are
\begin{equation}\label{nepo}
\dot p_j=-\partial_{q_j}H=-\left(\omega_j
+\frac{j^2}{2\pi\omega_j}
\Big(\sum_{\ell\in \N} \frac{\ell^2}{\omega_\ell}q_\ell^2\Big)\right)q_j\,,\qquad
\dot q_j=\partial_{p_j}H=\omega_j p_j
\,.
\end{equation}
\begin{rem}\label{pollastri}
By \eqref{nepo} if the initial data in \eqref{ginger}
have Fourier support  only on a subset $\mathcal S
 \subset \mathbb N$, namely
 $u_0(x)=\sum_{j\in\mathcal S}u_{0j}\phi_j(x)$
 and
  $v_0(x)=\sum_{j\in\mathcal S}v_{0j}\phi_j(x)$,
  then the same holds for the solution:
$u(t,x)=\sum_{j\in\mathcal S}u_{j}(t)\phi_j(x)$. 
In particular if the initial data are trigonometric polynomials the same holds for the solution. 
This proves Remark \ref{finitodim}.
See also Remark \ref{pollastre} below.                                                                                                                                                                                                                     
\end{rem}
Passing to complex coordinates, $\mathrm i=\sqrt{-1}\in \mathbb C$
\begin{equation}\label{tommaso}
z_j=\frac1{\sqrt 2}(q_j+\mathrm ip_j)\,\,\,\,\,\,\,\,\,\,\,\,\,\,\bar z_j=\frac1{\sqrt 2} (q_j-\mathrm i p_j)
\end{equation}
we get
\begin{equation}\label{mela}
H(z)=\Lambda(z) + G(z)\quad
\mbox{with}\quad
\Lambda(z):=\sum_{j\in \N}\omega_j z_j \bar z_j
\quad
G(z):=\frac{1}{32\pi}\Big( \sum_{j\in \N} \frac{j^2}{\omega_j}(z_j+\bar z_j)^2\Big)^2\,.
\end{equation}
We will also use the following notation
\begin{equation}\label{coordinate+-}
z_j^+:=z_j \ \ \ \ \ \ \ \ \ z_j^-:=\bar z_j \ \ \ \ \ \ \ \ \ \sigma:=(\sigma_1, \sigma_2)\in\{-,+\}^2\end{equation}
such that
\begin{equation*}
\begin{array}{c}
\sigma_h\sigma_k=\left\{
\begin{array}{l}
+\,\,\,\, \mathrm{if}\,\,\, \sigma_h=\sigma_k\\
-\,\,\,\, \mathrm{if}\,\,\, \sigma_h\neq\sigma_k
\end{array}
\right.\ \ \ \ \ \ \ \ \ \ \ \
\sigma_h+\sigma_k=\left\{
\begin{array}{ll}
0\,\,\,\, &\mathrm{if}\,\,\, \sigma_h\neq\sigma_k\\
2\,\,\,\,& \mathrm{if}\,\,\, \sigma_h=\sigma_k=+\\
-2\,\,\,\, &\mathrm{if}\,\,\, \sigma_h=\sigma_k=-
\end{array}
\right.
\end{array}
\end{equation*}
Set
$$\mathcal J:=\left\{(i,j,\sigma,\sigma')\in\N\times \N \times \{+,-\}^2\times \{+,-\}^2\ \right\}\,.
$$
Then we write the function $G$ in \eqref{mela} as
\begin{equation}\label{crostata}
G(z)
=\sum_{(i,j,\sigma,\sigma')\in \mathcal J}G_{ij}z_i^{\sigma_1} z_i^{\sigma_1'} z_j^{\sigma_2} z_j^{\sigma_2'}\quad
\mbox{where}\quad
0<G_{ij}=G_{ji}:=\frac 1{32\pi}\frac{i^2 j^2}{\omega_i \omega_j}
\leq \frac 1{32\pi}\,.
\end{equation}
Define the following subsets of
$\mathcal J$:
\begin{eqnarray}\label{poirot}
\mathcal I
&:=&
\left\{ (i,j,\s,\s')  :\   i\neq j \implies \sigma_1\neq \sigma_1'\,,\sigma_2\neq\sigma_2'\ \ {\rm and}\ \  i= j \implies  \sigma_1+\sigma_1'+\sigma_2+\sigma_2'=0 \right\}
\nonumber
\\
\mathcal I^c
&:=&
\mathcal J\setminus \mathcal I\,.
\end{eqnarray}
Then the \textsl{integrable part of $G$ is} 
\begin{equation}\label{prosciutto}
\bar G:=\sum_{(i,j,\sigma,\s')\in \mathcal I}G_{ij}z_i^{\sigma_1} z_i^{\sigma_1'} z_j^{\sigma_2} z_j^{\sigma_2'}=\sum_{i,j\in \N}A_{ij}|z_i|^2 |z_j|^2
\end{equation}
where
\begin{equation}\label{aerosol}
A_{ij}:=\left\{
\begin{array}{ll}
\displaystyle \frac 1{8\pi}\frac{i^2 j^2}{\omega_i \omega_j}\ \ \ \ &\mathrm{if}\ \ \ i\neq j\\
\\
\displaystyle \frac 3{16\pi}\frac{i^4}{\omega_i^2}\ \ \ \ &\mathrm{if}\ \ \ i=j\,.
\end{array}
\right.
\end{equation}

The idea of the fourth order BNF 
is to construct a 
canonical transformation 
eliminating  all the monomials in $G$
that are not in $\bar G$ so that, in the new variables,
the fourth order term reduces to $\bar G$ only.
The tool we will use here to construct 
canonical transformations is the \textsl{Lie series},
that we briefly discuss here \textsl{at a formal level}, 
the analytical estimates will be given in the next
section after having introduced the necessary
functional setting.
Define first the Poisson brackets  as
$$
\{F,G\}:=\im \sum_{j\in\N}
\big(\partial_{z_j}F\partial_{\bar z_j}G
-\partial_{\bar z_j}F\partial_{z_j}G
\big)\,.
$$
A typical way to produce symplectic transformation
is by the hamiltonian flow $\Phi^t_S$, for $t\in\R$, 
generated by an auxiliary Hamiltonian  $S$
(the ''generating function'').
Then the composition $H\circ\Phi^t_S$
can be written through  the Lie series
\begin{equation}\label{Lie}
H\circ\Phi^t_S=\sum_{n=0}^\infty
\frac{t^n}{n!}{\rm L}_S^n[H]\,,
\end{equation}
where  ${\rm L}_S[\cdot]:= \set{S,\cdot}$,
${\rm L}_S^0$ is the identity operator
and ${\rm L}_S^n:={\rm L}_S^{n-1}\circ{\rm L}_S$.
For the moment being the expression in \eqref{Lie}
is only formal, the convergence of the series
in suitable Hilbert spaces
will be discussed in Proposition
\ref{ham flow} below.

Coming back to
the fourth order BNF we have that if the 
auxiliary Hamiltonian $S$ satisfies
the ''homological equation'':
\begin{equation}\label{paprika}
G+\{S,\Lambda\}=\bar G\,,
\end{equation}
then by
\eqref{mela} and \eqref{Lie} (with $t=1$) we get
\begin{equation}\label{olio}
H\circ\Phi_S^1
=
\Lambda\circ\Phi_S^1+
G\circ\Phi_S^1=
\Lambda +G +\{S,\Lambda\}
+R
=
\Lambda +\bar G
+R\,,
\end{equation}
where $R$ contains only monomials of degree
greater or equal than 6.
Let us consider equation \eqref{paprika}.
 Noting that
\begin{equation}\label{facioli}
\left\{z_i^{\sigma_1}z_i^{\sigma_1'},z^+_i z_i^-\right\}=\mathrm i (\sigma_1+\sigma_1')z_i^{\sigma_1}z_i^{\sigma_1'}
\end{equation}
we have that
\begin{equation}\label{nuggets}
\left\{z_i^{\sigma_1}z_i^{\sigma_1'}z_j^{\sigma_2}z_j^{\sigma_2'},\Lambda\right\}
=\mathrm i\Delta_{ij}^{\sigma,\sigma'}
z_i^{\sigma_1}z_i^{\sigma_1'}z_j^{\sigma_2}z_j^{\sigma_2'}
\end{equation}
where 
$$
\Delta_{ij}^{\sigma,\sigma'}:=(\sigma_1+\sigma_1')\omega_i+(\sigma_2+\sigma_2')\omega_j\,.
$$
Set
\begin{equation}\label{linguine}
S:=\sum_{(i,j,\sigma,\sigma')\in \mathcal I^c}S_{ij}^{\sigma,\sigma'} z_i^{\sigma_1}z_i^{\sigma_1'}z_j^{\sigma_2}z_j^{\sigma_2'}\,,
\,\,\,\,\,\,\,\,\,\,\,\,
S_{ij}^{\sigma,\sigma'}:=\frac{\im G_{ij}}{\Delta_{ij}^{\sigma,\sigma'}}\,,
\end{equation}
where $\mathcal I^c$ was defined in \eqref{poirot}.
By \eqref{prosciutto}, \eqref{nuggets} and \eqref{linguine}
the auxiliary Hamiltonian $S$ solves 
\eqref{paprika} at formal level.

\section{Resonances}\label{sec3}

\subsection{Fourth order resonances}

To obtain the convergence of the series in 
\eqref{linguine} we have to 
estimates from below the fourth order
resonances  $\Delta_{ij}^{\sigma,\sigma'}$.

Recalling that we are assuming $m\geq -1/2$, we have that
\begin{equation}\label{chopin}
i>j\geq 1	\ \ 
\implies \ \ 
\o_i>\o_j\,.
\end{equation}

\begin{pro}\label{pera} 
We have\footnote{$\mu$ was defined in \eqref{trescul0}.}
\begin{equation}\label{cipster}
|\Delta_{ij}^{\sigma,\sigma'}|\geq \mu\,, \qquad
\forall\, 
(i,j,\sigma,\sigma')\in \mathcal I^c\,.
\end{equation}
\end{pro}
\noindent
\proof
We consider only the case $i\geq j$, the case
$i<j$ being analogous.
Then $\o_i\geq \o_j$.
Recalling \eqref{poirot}, for 
$(i,j,\sigma,\sigma')\in \mathcal I^c$
we have three cases:
\\
1) $\sigma_1=\sigma_1'$, $\sigma_2\neq\sigma_2'$ or
$\sigma_1\neq\sigma_1'$, $\sigma_2=\sigma_2'$\;
\\
2) $\sigma_1=\sigma_1'$, $\sigma_2=\sigma_2'$, $\sigma_1\neq\sigma_2$,  $i>j$
\\
3) $\sigma_1=\sigma_1'=\sigma_2=\sigma_2'$.
\\
Case 1). We have
$$
|\Delta_{ij}^{\sigma,\sigma'}|\geq 2 \o_j\geq 2\sqrt{1+m}
\geq 	\mu\,.
$$
Case 2). We have
\begin{eqnarray*}
|\Delta_{ij}^{\sigma,\sigma'}|
&=&
2(\omega_i-\omega_j)
=2\frac{\omega_i^2-\omega_j^2}{\omega_i+\omega_j}\geq \frac{\omega_i^2-\omega_j^2}{\omega_i}
=\frac{i^4-j^4+m(i^2-j^2)}{\omega_i}
\\
&=&
\frac{(i^2-j^2)(i^2+j^2+m)}{\o_i}
\geq 
\frac{(i^2-j^2)\o_i^2}{\o_i i^2}=
\frac{(i^2-j^2)\o_i}{ i^2}
=
\frac{(i^2-j^2)\sqrt{i^2+m}}{ i}
\,.
\end{eqnarray*}
 Since $i>j\geq 1$,  then $i\geq 2$
 and 
$$
|\Delta_{ij}^{\sigma,\sigma'}|
\geq 
\frac{(i^2-j^2)\sqrt{4+m}}{i}\geq \frac{[i^2-(i-1)^2]\sqrt{4+m}}{i}
=\frac {2i-1}i\sqrt{4+m}\geq 
\frac{3\sqrt{4+m}}{2}\geq  \mu\,.
$$
Case 3). We finally have
$$
|\Delta_{ij}^{\sigma,\sigma'}|=2(\omega_i+\omega_j)\geq 2\o_i\geq 2\sqrt{1+m}\geq \mu\,,
$$
concluding the proof.
\eproof

\subsection{Sixth order resonances}

 Recalling \eqref{facioli}
we have that
\begin{equation}\label{nuggets2}
\left\{z_i^{\sigma_1}z_i^{\sigma_1'}z_j^{\sigma_2}z_j^{\sigma_2'}
z_k^{\sigma_3}z_k^{\sigma_3'},\Lambda\right\}
=\mathrm i\Delta_{ijk}^{\sigma,\sigma'}
z_i^{\sigma_1}z_i^{\sigma_1'}z_j^{\sigma_2}z_j^{\sigma_2'}
z_k^{\sigma_3}z_k^{\sigma_3'}
\end{equation}
where 
$$
\Delta_{ijk}^{\sigma,\sigma'}:=(\sigma_1+\sigma_1')\omega_i+(\sigma_2+\sigma_2')\omega_j
+
(\sigma_3+\sigma_3')\omega_k\,.
$$
Let us introduce the following complementary
subsets of 
$
\mathbb N\times\mathbb N\times\mathbb N\times \{+,-\}^3\times \{+,-\}^3
$:
\begin{equation}\label{caciotta}
\Upsilon:=
\left\{(i,j,k,\s,\s') \ \ :\ \ \Delta_{ijk}^{\sigma,\sigma'}=0
\right\}\,,\qquad
\Upsilon^c:=
\left\{(i,j,k,\s,\s') \ \ :\ \ \Delta_{ijk}^{\sigma,\sigma'}\neq0
\right\}\,.
\end{equation}

\noi
The integrable part of the sixth order term of $R$ in \eqref{olio} is
$$
\bar R(z,\bar z):=
\sum_{(i,j,k,\sigma,\sigma')\in \Upsilon}R_{ijk}^{\sigma,\sigma'} z_i^{\sigma_1}z_i^{\sigma_1'}z_j^{\sigma_2}z_j^{\sigma_2'}
z_k^{\sigma_3}z_k^{\sigma_3'}
=\sum_{i,j,k}A_{ijk}|z_i|^2 |z_j|^2 |z_k|^2\,,
$$
for suitable $A_{ijk}\in \mathbb R$.
Set
\begin{equation}\label{linguine2}
\mathcal S:=\sum_{(i,j,k,\sigma,\sigma')\in \Upsilon^c}
\mathcal S_{ijk}^{\sigma,\sigma'} z_i^{\sigma_1}z_i^{\sigma_1'}z_j^{\sigma_2}z_j^{\sigma_2'}
z_k^{\sigma_3}z_k^{\sigma_3'}\,,
\,\,\,\,\,\,\,\,\,\,\,\,
\mathcal S_{ijk}^{\sigma,\sigma'}:=
\frac{\im R_{ijk}^{\sigma,\sigma'}}{\Delta_{ijk}^{\sigma,\sigma'}}\,.
\end{equation}
Now the homological equation is
\begin{equation}\label{paprika2}
R^{(6)}+\{\mathcal S,\Lambda\}=\bar R\,.
\end{equation}

\begin{pro}\label{pera2} 
Assume that $-1/2\leq m\leq 1$.
Then
\begin{equation}\label{cipster2}
|\Delta_{ijk}^{\sigma,\sigma'}|\geq 
 \frac78 |m|
\, \qquad
\forall\, 
(i,j,k,\sigma,\sigma')\in \Upsilon^c\,.
\end{equation}
\end{pro}
\proof
First note that 
\begin{equation}\label{vomito}
i\geq 1\,, \ -1/2\leq m\leq 1\quad\Longrightarrow\quad
\o_i\geq 1/\sqrt 2\,.
\end{equation}

Since  $(i,j,k,\sigma,\sigma')\in\Upsilon^c$
we are in one of the following three cases:
\begin{enumerate}
\item $|\Delta_{ijk}^{\sigma,\sigma'}|
=2|\omega_i+ \omega_j+ \omega_k|$, with $i\leq j\leq k$;

\item $|\Delta_{ijk}^{\sigma,\sigma'}|
=2|\omega_i- \omega_j+ \omega_k|$, with $i< j< k$;

\item $|\Delta_{ijk}^{\sigma,\sigma'}|
=2|\omega_i+ \omega_j- \omega_k|$, with $i\leq j< k$;

\item $|\Delta_{ijk}^{\sigma,\sigma'}|
=2|\omega_i- \omega_j- \omega_k|$, with $i< j\leq k$;

\item $|\Delta_{ijk}^{\sigma,\sigma'}|
=2|\omega_i+ \omega_j|$, with $i\leq j$;

\item $|\Delta_{ijk}^{\sigma,\sigma'}|
=2|\omega_i- \omega_j|$, with $i< j$;

\item $|\Delta_{ijk}^{\sigma,\sigma'}|
=2|\omega_i|$;

\end{enumerate}
In the cases 1 and 4  we have
$|\Delta_{ijk}^{\sigma,\sigma'}|\geq 2|\omega_k|
\geq \sqrt 2$ by \eqref{vomito}.
In the case 2 
we have
$|\Delta_{ijk}^{\sigma,\sigma'}|\geq 2|\omega_i|
\geq \sqrt 2$ again by \eqref{vomito}.
The cases 5 and 7 follow analogously.
The case 6 follows by Proposition \ref{pera}
noting that $\mu$ in \eqref{trescul0} satisfies
$$
\mu \geq \frac78 |m|
$$
for every $-1/2\leq m\leq 1$.
The only non trivial case is 3.
If $i=1$, since $k\geq j+1$ 
by \eqref{chopin} we have\footnote{Since
the minimum of $\omega_{j+1}-\omega_j-\sqrt{1+m}$ for $j\geq 1$ and 
$ m\geq -1/2$ is obtained for $j=1$
and $m=-1/2$.}
$$
\omega_k-\omega_j-\omega_i
\geq
\omega_{j+1}-\omega_j-\sqrt{1+m}
\geq \sqrt{14}-\sqrt 2
> 2.
$$
Assume now that $2\leq i\leq j$.
We distinguish two cases: 
$0\leq m\leq 1$ and $-1/2\leq m<0$.

Let us consider first the case 
$0\leq m\leq 1$.
Note that
\begin{equation}\label{mate}
-\frac{x^2}8\leq \sqrt{1+x}-1-\frac{x}2\leq 0\qquad \forall x\geq 0\,.
\end{equation}
Then for every $j\geq 1$ there exists a 
suitable $0\leq r_j\leq 1/8$ such that
$$
\omega_j=j^2\sqrt{1+\frac {m}{j^2}}
=j^2\left(1+\frac m{2j^2}-r_j \frac{m^2}{j^4} \right)=j^2+\frac m2 - r_j \frac{m^2}{j^2}\,.
$$
Analogously for $i$ and $k$.
We have
\begin{equation}\label{sonno}
\omega_i+\omega_j-\omega_k=i^2+j^2-k^2+\frac m2-r_i\frac{m^2}{i^2}-r_j\frac{m^2}{j^2}+r_k\frac{m^2}{k^2}
\end{equation}
Since $i^2+j^2-k^2\in \mathbb Z$
we distinguish two cases:
$i^2+j^2-k^2\leq -1$ and 
$i^2+j^2-k^2\geq 0$.
In the first case\footnote{Note that $k\geq 3$.}
$$
|\omega_i+\omega_j-\omega_k|
\geq 
1-\frac m2-\frac{m^2}{8k^2}
\geq
\frac12-\frac1{72}=\frac{35}{72}.
$$
In the second case, recalling that
$2\leq i\leq j$, we get
$$
\omega_i+\omega_j-\omega_k\geq
\frac m2-r_i\frac{m^2}{i^2}-r_j\frac{m^2}{j^2}
\geq m\left(
\frac12-\frac1{16}
\right)=\frac{7}{16}m\,,
$$
concluding the proof in the case $m$ positive.

Let us consider now the case
$-1/2\leq m<0$.
Note that
\begin{equation}\label{mate2}
-(3-2\sqrt 2)x^2\leq \sqrt{1+x}-1-\frac{x}2\leq 0\qquad \forall -1/2\leq x<0\,.
\end{equation}
Then for every $j\geq 1$ there exists a 
suitable $0\leq r_j'\leq 3-2\sqrt 2$ such that
$$
\omega_j=j^2\sqrt{1+\frac {m}{j^2}}
=j^2\left(1+\frac m{2j^2}-r_j' \frac{m^2}{j^4} \right)=j^2+\frac m2 - r_j' \frac{m^2}{j^2}\,.
$$
Analogously for $i$ and $k$.
We have
\begin{equation}\label{sonno2}
\omega_i+\omega_j-\omega_k=i^2+j^2-k^2+\frac m2-r_i'\frac{m^2}{i^2}-r_j'\frac{m^2}{j^2}+r_k'\frac{m^2}{k^2}
\end{equation}
Since $i^2+j^2-k^2\in \mathbb Z$
we distinguish two cases:
$i^2+j^2-k^2\leq 0$ and
$i^2+j^2-k^2\geq 1$.
In the first case, recalling that $k\geq 3$,
$$
|\omega_i+\omega_j-\omega_k|
\geq 
-\frac m2\left(1-\frac{3-2\sqrt 2}{9}
\right)
\geq
-0.49\, m.
$$
In the second case, recalling that
$2\leq i\leq j$, we get by \eqref{sonno2}
$$
\omega_i+\omega_j-\omega_k\geq
1+
\frac m2-r_i'\frac{m^2}{i^2}-r_j'\frac{m^2}{j^2}
\geq 
1-\frac14-\frac18=\frac58
$$
concluding the proof also 
in the case $m$ negative.
\eproof


\section{Quantitative Birkhoff Normal Forms}\label{sec4}

\subsection{Weighted Sobolev norms and
improved stability results}

In order to optimize the estimates we introduce equivalent norms on $H^s_*$.
Given $N\geq 1$ let us define the two norms
\begin{equation}\label{enorme}
|u|_{H^{s}_N}^2=
\sum_{j\leq N}j^{2}u_{j}^2
+\sum_{j> N}
j^{2s}u_{j}^2
\,,\quad 
|u|_{\tilde H^{s}_N}^2
=
\sum_{j\leq N}\frac{j^{2}}{\o_j^2}u_{j}^2
+\sum_{j> N}
\frac{j^{2s+4}}{\o_j^2}u_{j}^2\,,
\end{equation}
which are equivalent 
to the $H^s$--norm; indeed for every $u\in H^s_*$, $s\geq 1$, one has
$$
N^{1-s}|u|_s
 \leq
|u|_{H^{s}_N}
\leq 
|u|_s
$$
and, for $s>0$,
$$
\min\{\o_N^{-1}N^{1-s},\o_N^{-1} N^2,1\} |u|_s
\leq
|u|_{\tilde H^{s}_N}
\leq 
|u|_s\,.
$$
Note that if 
$u$ is a trigonometric polynomial
of degree $N$, then  $|u|_{H^{s}_N}=|u|_1$.
In general, if $u\in H^s_*$ for every $\delta>0$ there exists $N$
(depending on $\delta,s$ and $u$) such that
\begin{equation}\label{CEV}
|u|_{H^s_N}\leq (1+\delta)|u|_1
\end{equation}
and, analogously,
\begin{equation}\label{CEV2}
|u|_{\tilde H^s_N}\leq (1+\delta)|u|_{\tilde H^{-1}}\,,\ \ \ \mbox{with}\ \ \ 
|u|_{\tilde H^{-1}}^2:=\sum_{j\geq 1}\frac{j^{2}}{\o_j^2}u_{j}^2\,.
\end{equation}
Finally if $s\geq 1$, by \eqref{mac}
we get
\begin{equation}\label{geppetto}
\sup_{x\in[0,\pi]}|u(x)|\leq c_1|u|_1
=\sqrt{\pi/3}|u|_1
\leq 
\sqrt{\pi/3}|u|_{H^1_N}
\leq
\sqrt{\pi/3}|u|_{H^s_N}\,.
\end{equation}
By \eqref{CEV}-\eqref{geppetto}
it is immediate to see that the following result is a stronger formulation of 
Theorem \ref{patata}:
\begin{thm}\label{patata2}
Let $m\geq -1/2$, $s>1/2$ and assume that
the initial data in \eqref{pisolino}
satisfy
$u_0\in H^{s+2}_*$, 
$v_0\in H^{s}_*$ 
and
 \begin{equation}\label{trescul03}
\e:=\max\{|u_0|_1\,,\ 
|v_0|_{\tilde H^{-1}}\}\ \leq \e_0:=
0.08\sqrt\mu
 \,.
\end{equation}
Then there exists $N$ such that 
the weak solution $u(t,x)$
of equation \eqref{checco}-\eqref{ginger}
satisfies
\begin{equation}
u(t,\cdot)\in H^{s+2}_*\,,\ \ 
|u(t,\cdot)|_{H^{s+2}_N}\leq
  1.9 \e\,,\ \ 
 \sup_{x\in[0,\pi]}|u(t,x)| 
  \leq 2\e\,,\ \ 
\forall\, |t|\leq T_0\mu\e^{-4}\,,
\label{burrata}
\end{equation}
where
\begin{equation}\label{stoccafisso}
T_0=0.948
\frac{1-c_\dag c_*\e^2/\mu}{c_\dag c_*^2}
\,,\quad
c_*:=
\left(\frac{11 }{10}\right)^2
\left(\frac56-\frac{1}{80e}\right)^{-2}
\,,\qquad
c_\dag:=\frac{44e}{\pi}\,.
\end{equation}
\end{thm}
\begin{rem}
(i) Note that $\e$ defined in \eqref{trescul03}
 is lesser or equal than $\e$ defined in \eqref{trescul0}
 and \eqref{burrata} implies \eqref{burrata0}.
 In the following we will use the definition of
 $\e$ in \eqref{trescul03} and prove Theorem \ref{patata2}
 from which Theorem \ref{patata} follows.
 \\
 (ii) The choice of $N$ and, therefore,
 the one of the norm 
 $|\cdot|_{H^{s+2}_N}$ (which is equivalent to the standard
 norm $|\cdot|_{H^{s+2}}$)
  depend on the initial 
 data $u_0$ and $v_0$.
\end{rem}

Analogously, recalling the definitions of 
$\e_1$ and $T_1$ in \eqref{e1},
 we give a stronger formulation of Theorem \ref{patataBIS}:

\begin{thm}\label{patataTER}
In addition to the assumptions of Theorem
\ref{patata2}
 we require that
 $-1/2\leq m\leq 1$, $m\neq 0$ and $\e\leq \e_1$.
Then
there exists $N$ such that 
\begin{equation}
u(t,\cdot)\in H^{s+2}_*\,,\ \ 
|u(t,\cdot)|_{H^{s+2}_N}\leq
  2 \e\,,\ \ 
 \sup_{x\in[0,\pi]}|u(t,x)| 
  \leq 2.1\e\,,\ \ 
\forall\, |t|\leq T_1\e^{-6}\,.
\label{burrataBIS}
\end{equation}
\end{thm}


We assume that the coordinates $z$ introduced in \eqref{tommaso} belong to some complete subspace  of $\ell^2$.
More precisely, given a real sequence $\mathtt w=\{\mathtt w_j\}_{j\in\N},$
with $\inf_{j\in\N}\mathtt w_j>0$, we consider the 
Hilbert space\footnote{
	Endowed with the scalar product
	$(z,v)_{\th_\tw}:=\sum_{j\in\N} \mathtt w_j^2 z_j \bar v_j.$}
\begin{equation}\label{pistacchio}
\th_{\mathtt w}
:=
\set{z:= \pa{z_j}_{j\in\N}\in\ell^2(\C)\,: \quad \abs{z}_{\tw}^2
	:= 
	\sum_{j\in\N} \mathtt w_j^2 \abs{z_j}^2 < \infty}
\end{equation}
and fix  the symplectic structure to be 
$
\im\sum_j d z_j\wedge d \bar z_j
$.
For example taking $\tw_j=\jap{j}^{s}$  we have $\th_{\mathtt w}=H^s(\mathbb C)$; 
then
by \eqref{nduja}-\eqref{tommaso} and since
$\o_j\sim j^2$,
\begin{equation}\label{peppa}
u(x)=\sum_{j\in \N}\frac{z_j+\bar z_j}{\sqrt{\pi\omega_j}}\sin(jx)
\in  H^{s+1}\,.
\end{equation} 
\smallskip

\noindent Here and in the following, given $r>0$, by $B_r(\th_{\mathtt w})$ we mean the closed
ball of radius $r$ centered at the origin of $\th_{\mathtt w}.$

\begin{defn}[Regular Hamiltonians]\label{Hr}
	For $ r>0$,  let
	$\mathcal{A}_r(\th_{\mathtt w})$
	be the space of  
	Hamiltonians 
	$$
	H : B_r(\th_{\mathtt w}) \to \R
	$$ 
	such that there exists a point-wise  absolutely convergent power series 
	expansion\footnote{As usual  
	$\a:=(\a_1,\a_2,\ldots)$ is a multi-index with $\a_j\in\N_0:=\{0,1,2,\ldots\}$
	and
		$|k|:=\sum_{j\in\N}|k_j|$.}
	\begin{equation}\label{fragola}
	H(z)  = \sum_{\substack{\a,\b\in\N_0^\N\,, \\
			|\a|+|\b|<\infty} }H_{\a,\b}z^\a \bar z^\b\,,
	\qquad
	z^\a:=\prod_{j\in\Z}u_j^{\a_j}
	\end{equation}
	with\footnote{This means that $H$ is real analytic in the real and imaginary part of $z$.} 
	$H_{\a,\b}= \overline{ H}_{\b,\a}$.
	We say that
	$H\in \mathcal{A}_r(\th_{\mathtt w})$  is  a majorant analytic Hamiltonian 
			if its majorant, namely
	$\und H:  B_r(\th_{\mathtt w}) \to \R$  as
	\begin{equation}\label{betta}
	\und H(z)  = \sum_{
			|\a|+|\b|<\infty}|H_{\a,\b}|z^\a \bar z^\b\,,
	\end{equation}
	is point-wise  absolutely convergent.
We denote by
$\cH_{r}(\th_{\mathtt w})$ 
the subspace
of all $H\in\mathcal{A}_r(\th_{\mathtt w})$
with $H(0)=0$
such that  
\begin{equation}\label{salame}
|H|_{r,\tw}
:=
r^{-1} \left(\sup_{\norm{u}_{\tw}\leq r} 
\norm{{X}_{\underline H}}_{\tw} \right) < \infty\,,
\end{equation}
where
${X}_{\underline H}$ is the 
(complex) Hamiltonian vector field of $\underline H$
whose entries are defined as
$
X_{\underline H}^{(j)}  = 
\im \frac{\partial}{\partial \bar z_j} \underline H
$ for $j\in\N$. $\cH_{r}(\th_{\mathtt w})$
endowed with the norm $|\cdot|_{r,\tw}$ is a Banach space.
\end{defn}
By \eqref{betta} and \eqref{salame}
it is simple to see that
$H$ and $H'$,  written as in \eqref{fragola},
satisfy the monotonicity property, namely                                         
\begin{equation}\label{monopoli}
|H_{\a,\b}|\leq |H'_{\a,\b}|\,,\ \forall\, \a,\b\ \ \Longrightarrow\ \ 
|H|_{r,\tw}\leq |H'|_{r,\tw}\,.
\end{equation}

We now show that $G$ defined in \eqref{mela}
is a regular Hamiltonian.
First we write it as in \eqref{fragola}, namely
$$
G(z)=\sum_{|\a|+|\b|=4}
G_{\a,\b}z^\a \bar z^\b\,,\quad
\mbox{for suitable}
\quad
0\leq G_{\a,\b}\leq
\frac{1}{4\pi}
\,. 
$$
For example\footnote{As usual
$e_i$ is the infinite vector whose 
$i$--th entry is 1 and all the others are 0.
}
$G_{e_i+e_j,e_i+e_j}$
is equal to $2 A_{ij}$ defined in \eqref{aerosol}
when $i\neq j$ and to $A_{ij}$ when $i=j$.
Note that, recalling \eqref{betta} and
since the coefficients $G_{\a,\b}$
are non negative,
$G$ coincides with its majorant, namely
 $G=\underline G$.
 Then, by \eqref{mela}, for $h\in\N$,
\begin{equation}\label{temperatura}
X_{\underline G}^{(h)} (z) = 
\im \frac{\partial}{\partial \bar z_h} G(z)=
\frac{\im}{8\pi}
\frac{h^2}{\omega_h}
(z_h+\bar z_h)
 \sum_{j\in \N} \frac{j^2}{\omega_j}(z_j+\bar z_j)^2
 \,.
\end{equation}
 Recalling \eqref{pistacchio}, let us consider now the sequence 
 \begin{equation}\label{peso}
\tw^0:=\{j/\sqrt{\o_j}\}_{j\in\N}\quad
\mbox{inducing the norm}\quad
|z|_{\tw^0}^2= \sum_{j\in \N} \frac{j^2}{\omega_j}|z_j|^2\,,
\end{equation}
which is equivalent to the $\ell^2$--norm 
$|z|_{\ell^2}^2=\sum_{j\in \N} |z_j|^2$,
since, recalling \eqref{e1},
\begin{equation}\label{golfetta}
0<\inf_{j\in\N}j^2/\o_j
\leq \sup_{j\in\N}j^2/\o_j=\cccp
\end{equation}
 Note that with $u$ as in \eqref{peppa}
 $$ \sum_{j\in \N} \frac{j^2}{\omega_j}(z_j+\bar z_j)^2
 =|z+\bar z|_{\tw^0}^2
 =|u_x|_{L^2}^2\,.
 $$
 Then for every weight $\tw$  we have that
by \eqref{temperatura} and \eqref{golfetta}
we get
$$
 |X_{\underline G}^{(h)} (z)|_{\tw}\leq \frac{\cccp}{8\pi}
 |z+\bar z|_{\tw}|z+\bar z|_{\tw^0}^2
\,.
$$
Recalling the definition in \eqref{salame}, for $r>0$ we have
\begin{equation}\label{mortadella}
 |G|_{r,\tw}\leq \frac{\cccp}{\pi}
 \g_r^2\,,\qquad \g_r:=\sup_{|z|_\tw\leq r}|z|_{\tw_0}\leq r\,.
\end{equation}
Note also that for every $\l>0$
\begin{equation}\label{tachi}
\g_{\l r}=
\l \g_r\,.
\end{equation}
Analogously by \eqref{monopoli}, \eqref{prosciutto} and \eqref{mortadella}
we get
\begin{equation}\label{mortadella2}
 |\bar G|_{r,\tw}\leq \frac{\cccp}{\pi}
 \g_r^2\,,
 \qquad
 |\bar G-G|_{r,\tw}\leq \frac{\cccp}{\pi}
 \g_r^2\,.
\end{equation}

\subsection{Estimates on 
Poisson brackets and Lie series}

The following result is proved in \cite{BMP}
(Proposition 2.1).

\begin{pro}\label{piadina}
	For $0 <\rho\leq r$  we have
	\begin{equation}\label{commXHK}
	|\{F,G\}|_{r,\tw}
	\le 
	4\pa{1+\frac{r}{\rho}}
	|F|_{r+\rho,\tw}
	|G|_{r+\rho,\tw}\,.
	\end{equation}
\end{pro}

The following result is a consequence
 of Proposition
\ref{piadina}.
It is essentially  Lemma 2.1 of \cite{BMP} 
with improved constants. 
Its proof is given in the Appendix.
\begin{pro}\label{ham flow}
	Let $0<\rho< r $,  and $S\in\cH_{r+\rho}(\th_{\mathtt w})$ with 
	\begin{equation}\label{stima generatrice}
	4e \left(1+\frac{r}{\rho}\right)\abs{S}_{r+\rho,\tw} :=
	\eta<1\,. 
	\end{equation} 
	Then the time $1$-Hamiltonian flow 
	$\Phi^1_S: B_r(\th_{\mathtt w})\to
	B_{r + \rho}(\th_{\mathtt w})$  is well defined, analytic, symplectic with
	\begin{equation}
	\label{pollon}
	\sup_{u\in  B_r(\th_{\mathtt w})} 	\norm{\Phi^1_S(u)-u}_{\th_{\mathtt w}}
	\le
	(r+\rho)  \abs{S}_{r+\rho,\tw}
	\leq
	\frac{\rho\eta}{4 e}.
	\end{equation}
	For any $H\in \cH_{r+\rho}(\th_{\mathtt w})$
	we have that
	$H\circ\Phi^1_S\in\cH_{r}(\th_{\mathtt w})$ and
	\begin{equation}
	\label{caio}
	\abs{H\circ\Phi^1_S-H}_{r,\tw}
	\le \frac{1}{e} \frac{\eta}{1-\eta}
	\abs{H}_{r+\rho,\tw}\,.
	\end{equation}
	In general, for any 
	 sequence  $\{c_n\}_{n\in\N}$, we have 
	\begin{equation}\label{brubeck}
	\Big|\sum_{n\geq h} c_n {\rm L}^n_S[H]\Big|_{r,\tw} \le
	\sup_{n\geq h}  \big(|c_n|n^n e^{-n}\big)
	\frac{\eta^h}{1-\eta}
	 |H|_{r+\rho,\tw} 
	\,.
	\end{equation}
\end{pro}


\subsection{Fourth and sixth order BNF}

By Proposition \ref{pera}, \eqref{linguine}, \eqref{monopoli} and \eqref{mortadella}
we get
\begin{equation}\label{mortadella3}
 |S|_{r,\tw}\leq \frac{\cccp}{\pi \mu}
 \g_r^2\,.
\end{equation}
\begin{thm}[$4^{th}$--order Birkhoff Normal Form]\label{sob}
	Let $r>0$ such that $\g_r$ defined in \eqref{mortadella}
	satisfies
	\begin{equation}\label{daje}
0<\g_r\leq \frac{5}{11}\sqrt{\frac{\pi\mu}{22e\cccp}}\,.
\end{equation}
	 Then
	there exist two close--to--the--identity invertible symplectic transformations 
	$$
	\Psi,\Psi^{-1}:\quad B_{r}(\th_{\tw})
	\mapsto \th_{\tw}
	$$
	satisfying
	\begin{eqnarray}\label{stracchino}
	&   
	\sup_{|z|_{\tw}\leq r}
	|\Psi^{\pm 1}(z)-z|_{\tw} \le 
	\frac{11^3\cccp}{10^3\pi\mu}
	 \g_r^2 \,r 
	\leq \frac{1}{80e} r
	\\
	&\Psi\big(\Psi^{-1}(z)\big)= \Psi^{-1}\big(\Psi(z)\big)= z \,,\quad \forall z\in B_{(1-\frac{1}{80e}) r}(\th_{\tw})
	\label{stracchinobis}
	\end{eqnarray}
	such that in the new coordinates
	\begin{equation}\label{pesce}
	\mathrm H:=
	H\circ \Psi= \Lambda + \bar G + R
	\end{equation}
	with
	\begin{equation}\label{duspaghibis}
	|R|_{r,\tw}\leq C_r:=
	\frac{3\cccp}{2\pi e}\frac{\eta_r}{1-\eta_r}
 \left(\frac{11 \g_r}{10}\right)^2\,, \quad
 \mbox{with}\quad
 \eta_r:=
\frac{44e\cccp}{\pi\mu}
\left(\frac{11 \g_r}{10}\right)^2	
\end{equation}
Moreover $R=\sum_{n\geq 3} R^{(2n)}$ 
is a series of homogeneous
polynomials of even degree $2n$ 
of the form
\begin{equation}\label{cacao}
R^{(2n)}(z,\bar z)=\sum_{\ell\in \mathbb N^n\!\!,\, \,\sigma, \sigma'\in \{+,-\}^n}
R^{\sigma,\sigma'}_\ell
z_{\ell_1}^{\sigma_1}z_{\ell_1}^{\sigma_1'}
\cdots
z_{\ell_n}^{\sigma_n}z_{\ell_n}^{\sigma_n'}
\end{equation}
and, setting $R^{(\geq 8)}:=R-R^{(6)}$, we have
\begin{equation}\label{baccala}
|R^{(\geq 8)}|_{r,\tw}\leq \frac{16\cccp}{3e^2\pi^3}
\frac{\g_r^6}{\mu^2} 
\end{equation}
\end{thm}
\begin{rem}\label{pollastre}
(i)
In view of \eqref{cacao} the Hamilton equations of 
$\mathrm H$ in \eqref{pesce} have the same structure as
\eqref{nepo}, namely the terms $z_j^\pm$ factorize in $\partial_{z_j}\mathrm H$.
Therefore, as for $H$,
 if the initial data 
have Fourier support  only on a subset $\mathcal S
 \subset \mathbb N$,
  then the same holds for the solution.   
  The same statement will hold for the Hamiltonian $\mathtt H$ in \eqref{pesce6}
  in view of \eqref{cacao6}.    
  \\
  (ii) Partial BNF are given in 
\cite{BB05}
\cite{GY}) 
  for the beam equation but with  nonlinearities
different from the one in \eqref{checco}.                                                                                                                                                                                                             
\end{rem}
\proof
We use Proposition \ref{ham flow} with $\rho=r/10$.
Set $\Psi^{\pm 1}:=\Phi^{\pm 1}_S$ with $S$ defined in 
\eqref{linguine}.
By \eqref{mortadella3} (used with
$r\to 11r/10$ and recalling \eqref{tachi}) 
and \eqref{daje} we have that 
the smallness condition on $S$ contained in 
\eqref{stima generatrice} holds, indeed
\begin{equation}\label{bagno}
\eta=
44e |S|_{\frac{11}{10}r,\tw}\leq \eta_r=
44e\cccp\frac{11^2}{10^2\pi\mu}
 \g_r^2\leq \frac12\,.
\end{equation}
Then Proposition \ref{ham flow} applies.
By \eqref{pollon} and \eqref{bagno} we get
\eqref{stracchino} for $\Psi$.
An analogous estimate holds for $\Psi^{-1}$.
\eqref{stracchinobis} follows by \eqref{stracchino}.
\\
Now note that, by \eqref{paprika}, 
$\{S,\Lambda\}=\bar G-G$. Recalling 
\eqref{olio}
\begin{eqnarray*}
H\circ\Phi_S^1
&=&
\Lambda +G +\{S,\Lambda\} 
+\sum_{n=2}^\infty
\frac{1}{n!}{\rm L}_S^n[\Lambda]
+ (G\circ\Phi_S^1-G)
\\
&=&
\Lambda +\bar G
+\sum_{n=2}^\infty
\frac{1}{n!}{\rm L}_S^{n-1}[\bar G-G]
+ (G\circ\Phi_S^1-G)
\\
&=&
\Lambda +\bar G+R\,,
\end{eqnarray*}
where
\begin{equation}\label{bio}
R:=R_1+R_2\,,\quad
R_1:=\sum_{n=1}^\infty
\frac{1}{(n+1)!}{\rm L}_S^{n}[\bar G-G]
\,,
\quad
R_2:= (G\circ\Phi_S^1-G)\,.
\end{equation}
By \eqref{brubeck}, \eqref{bagno} and \eqref{mortadella2} (with $r\to r+\rho=11r/10$)
 we get
$$
|R_1|_{r,\tw}\leq \frac{1}{2e}\frac{\eta_r}{1-\eta_r}
\frac{\cccp}{\pi} \left(\frac{11 \g_r}{10}\right)^2\,.
$$
By \eqref{caio} we get
$$
|R_2|_{r,\tw}\leq 
\frac{1}{e}\frac{\eta_r}{1-\eta_r}
\frac{\cccp}{\pi} \left(\frac{11 \g_r}{10}\right)^2\,.
$$
Then \eqref{duspaghibis} follows.
Moreover \eqref{cacao} follows by induction 
and \eqref{Lie} (with $t=1$). 
Finally, noting that
$$
R^{(\geq 8)}:=
\sum_{n=2}^\infty
\frac{1}{(n+1)!}{\rm L}_S^{n}[\bar G-G]
+\sum_{n=2}^\infty\frac{1}{n!}{\rm L}_S^{n}[G]\,,
$$
 by
\eqref{brubeck}, \eqref{mortadella}, \eqref{mortadella2}, \eqref{mortadella3} and  we get \eqref{baccala}.
\eproof


By Proposition \ref{pera2}, \eqref{monopoli} and 
\eqref{duspaghibis}
we get
\begin{equation}\label{mortadella4}
 |\mathcal S|_{r,\tw}\leq \frac{8 C_r}{7|m|}
 \,.
\end{equation}

\begin{thm}[$6^{th}$--order Birkhoff Normal Form]\label{sob6}
	Let $\tr:=10 r/11$ such that\footnote{Recall \eqref{duspaghibis}.} 
		\begin{equation}\label{daje6}
 \tilde\eta_\tr:=
\frac{352 e C_{11\tr/10}}{7|m|}
 \leq \frac12\,.
\end{equation}
	 Then 
	$$
	\Phi_{\mathcal S}^1,\Phi_{\mathcal S}^{-1}:
	\quad B_{\tr}(\th_{\tw})
	\mapsto \th_{\tw}
	$$
	satisfying
	\begin{eqnarray}\label{stracchino6}
	&   
	\sup_{|z|_{\tw}\leq \tr}
	|\Phi_{\mathcal S}^{\pm 1}(z)-z|_{\tw} 
	\leq \frac{1}{80e} \tr
	\\
	&\Phi_{\mathcal S}^{1}\big(\Phi_{\mathcal S}^{-1}(z)\big)= \Phi_{\mathcal S}^{-1}
	\big(\Phi_{\mathcal S}^{1}(z)\big)= z \,,\quad \forall z\in B_{(1-\frac{1}{80e}) \tr}(\th_{\tw})
	\label{stracchinobis6}
	\end{eqnarray}
	such that in the new coordinates
	\begin{equation}\label{pesce6}
	\mathtt H:=
	\mathrm H\circ\Phi_{\mathcal S}^1= \Lambda + \bar G + \bar R
         +\mathtt R
	\end{equation}
	with
	\begin{equation}\label{duspaghibis6}
	|\mathtt R|_{\tr,\tw}\leq \tilde C_\tr:=
	\frac{1}{44e^2}\tilde\eta_\tr^2+	
	\frac{2}{e}\tilde\eta_\tr
\left(\frac{\cccp}\pi \g_{11\tr/10}^2
+C_{11\tr/10}
\right)+
\frac{16\cccp}{3e^2\pi^3}
\frac{\g_{11\tr/10}^6}{\mu^2} 
(1+2\tilde\eta_\tr/e)\,.
\end{equation}
Finally $\mathtt R=\sum_{n\geq 4} \mathtt R^{(2n)}$ 
is a series of homogeneous
polynomials of even degree $2n$ 
of the form
\begin{equation}\label{cacao6}
\mathtt R^{(2n)}(z,\bar z)=\sum_{\ell\in \mathbb N^n\!\!,\, \,\sigma, \sigma'\in \{+,-\}^n}
\mathtt R^{\sigma,\sigma'}_\ell
z_{\ell_1}^{\sigma_1}z_{\ell_1}^{\sigma_1'}
\cdots
z_{\ell_n}^{\sigma_n}z_{\ell_n}^{\sigma_n'}
\end{equation}
\end{thm}
\proof
We use Proposition \ref{ham flow} with
$S\to \mathcal S$,
 $r\to \tr$ and $\rho\to\tr/10$.
By \eqref{mortadella4}  
and \eqref{daje6} we have that 
the smallness condition on $\mathcal S$ contained in 
\eqref{stima generatrice} holds, indeed
\begin{equation}\label{bagno6}
\eta:=
44e |\mathcal S|_{\frac{11}{10}\tr,\tw}\leq \tilde\eta_\tr\,.
\end{equation}
Then Proposition \ref{ham flow} applies.
By \eqref{pollon} and \eqref{bagno6} we get
\eqref{stracchino6} for $\Phi_{\mathcal S}^{1}$.
An analogous estimate holds for $\Phi_{\mathcal S}^{-1}$.
\eqref{stracchinobis6} follows by \eqref{stracchino6}.
\\
Now note that, by \eqref{paprika2}, 
$\{\mathcal S,\Lambda\}=\bar R-R^{(6)}$.
 Recalling 
 \eqref{mela}, \eqref{pesce}, \eqref{baccala} and \eqref{Lie} with $t=1$, we get
\begin{eqnarray*}
&&\mathrm H\circ\Phi_{\mathcal S}^1
=
\Lambda\circ\Phi_{\mathcal S}^1+
\bar G\circ\Phi_{\mathcal S}^1+R^{(6)}\circ\Phi_{\mathcal S}^1
+R^{(\geq 8)}\circ\Phi_{\mathcal S}^1
\\
&=&
\Lambda +\bar G +\{{\mathcal S},\Lambda\} + R^{(6)}+
\sum_{n=2}^\infty
\frac{1}{n!}{\rm L}_{\mathcal S}^n[\Lambda]
+ \Big((\bar G+R^{(6)})\circ\Phi_{\mathcal S}^1-(\bar G+R^{(6)})\Big)
+R^{(\geq 8)}\circ\Phi_{\mathcal S}^1
\\
&=&
\Lambda +\bar G +\{{\mathcal S},\Lambda\} + R^{(6)}+
\sum_{n=2}^\infty
\frac{1}{n!}{\rm L}_{\mathcal S}^{n-1}[\bar R-R^{(6)}]
+ \Big((\bar G+R^{(6)})\circ\Phi_{\mathcal S}^1-(\bar G+R^{(6)})\Big)
+R^{(\geq 8)}\circ\Phi_{\mathcal S}^1
\\
&=&
\Lambda +\bar G+\bar R+\mathtt R\,,
\end{eqnarray*}
where
\begin{eqnarray*}
&&\mathtt R:=\mathtt R_1+\mathtt R_2+\mathtt R_3\,,\quad
\mathtt R_1:=\sum_{n=1}^\infty
\frac{1}{(n+1)!}{\rm L}_{\mathcal S}^{n}[\bar R-R^{(6)}]
\,,
\\
&&\mathtt R_2:= 
(\bar G+R^{(6)})\circ\Phi_{\mathcal S}^1-(\bar G+R^{(6)})
\,,
\qquad
\mathtt R_3:=R^{(\geq 8)}\circ\Phi_{\mathcal S}^1
\,.
\end{eqnarray*}
We first note that by \eqref{monopoli} and \eqref{duspaghibis}
$$
|\bar R-R^{(6)}|_{11\tr/10,\tw}\,,\ 
|R^{(6)}|_{11\tr/10,\tw}\,,\ 
|R^{(\geq 8)}|_{11\tr/10,\tw}\leq C_{11\tr/10}
$$
By \eqref{brubeck} and \eqref{bagno6} 
 we get
$$
|\mathtt R_1|_{\tr,\tw}\leq \frac{1}{44e^2}\tilde\eta_\tr^2
\,.
$$
By \eqref{caio} and \eqref{mortadella2} we get
$$
|\mathtt  R_2|_{\tr,\tw}\leq 
\frac{2}{e}\tilde\eta_\tr
\left(\frac{\cccp}\pi \g_{11\tr/10}^2
+C_{11\tr/10}
\right)\,.
$$

Writing
$\mathtt R_3:=R^{(\geq 8)}+(R^{(\geq 8)}\circ\Phi_{\mathcal S}^1-R^{(\geq 8)})$
 by  \eqref{caio} we get
$$
|\mathtt  R_3|_{\tr,\tw}\leq 
\frac{16\cccp}{3e^2\pi^3}
\frac{\g_{11\tr/10}^6}{\mu^2} 
(1+2\tilde\eta_\tr/e)
\leq
\frac{16(1+e)\cccp}{3e^3\pi^3}
\frac{\g_{11\tr/10}^6}{\mu^2} 
\,.
$$
Then \eqref{duspaghibis6} follows.
Finally \eqref{cacao6} follows by induction, \eqref{linguine2} 
and \eqref{Lie} (with $t=1$). 
\eproof

\section{Proof of the stability estimates}\label{sec5}

We need the following result, which is 
 Lemma 7.1 of \cite{BMP} with different constants.
For 
completeness we give the proof in the Appendix.

\begin{lem}\label{cobra}
 On the Hilbert space $\th_\tw$ consider the Hamilton equation
  \[
  \dot z = X_{\mathcal N}+X_R\,,
  \qquad
  z(0)=z_0\,, \qquad
  |z_0|_\tw< \frac56 r\,,
 \]
 where 
$\mathcal N\in \mathcal{A}_r(\th_{\mathtt w}) $
and $R\in\cH_{r}(\th_{\mathtt w})$
for some $r>0.$
Assume that 
$$
{\rm Re} (X_{\mathcal N},z)_{\th_\tw}=0\,.
$$  
Then
\begin{equation}\label{giuncata}
\Big| |z(t)|_\tw-|z_0|_\tw\Big|< \frac{r}{6}\,, \qquad
\forall\,  |t|\leq 
\frac{1}{6|R|_{r,\tw}}\,.
\end{equation}
\end{lem}


\bigskip

\noindent
\subsubsection*{Proof of Theorem \ref{patata2}}

We now apply the above lemma to 
the Hamiltonian $\mathrm H=\Lambda + \bar G + R$
in \eqref{pesce} with ${\mathcal N}=
\Lambda + \bar G$. Then, by \eqref{duspaghibis} we get the following
\begin{cor}\label{gina}
 Consider the Hamilton equation of the Hamiltonian ${\mathrm H}$
 in \eqref{pesce}: 
  \[
  \dot \z = X_{\mathrm H}\,,
  \qquad
  \z(0)=\z_0\,, \qquad
  |\z_0|_\tw<\frac56 r\,,
 \]
 with $r>0$
 such that $\g_r$ defined in \eqref{mortadella}
	satisfies \eqref{daje}.
 Then\footnote{Recall that
 $\eta_r$ is defined in \eqref{duspaghibis}.}
 \begin{equation}\label{giuncata2}
\Big| |\z(t)|_\tw-|\z_0|_\tw\Big|< \frac{r}{6}\,, 
\quad
|\z(t)|_\tw\leq r
\qquad
\forall\,  |t|\leq 
T:=\frac{\pi e}{9}\frac{1-\eta_r}{\eta_r}
 \left(\frac{10}{11 \g_r}\right)^2
\,.
\end{equation}
\end{cor}

\begin{cor}\label{moderno}
  Consider the Hamilton equation of the Hamiltonian $H$
 in \eqref{mela}: 
  \[
  \dot z = X_{H}\,,
  \qquad
  z(0)=z_0\,, \qquad
  |z_0|_\tw< \left(\frac56-\frac{1}{80e}\right) r\,,
 \]
 with $r>0$ 
 such that $\g_r$ defined in \eqref{mortadella}
	satisfies \eqref{daje}.
 Then
 \begin{equation}\label{giuncata3}
|z(t)|_\tw\leq
\left(1+\frac{1}{80e}\right) r
\qquad
\forall\,  |t|\leq 
 T
\,.
\end{equation}
\end{cor}
\proof
Recalling \eqref{stracchino}
we define $\z_0:=\Psi^{-1}(z_0)$
with 
$$
|\z_0|_\tw<\frac56 r\,.
$$
Then by \eqref{giuncata2}
and \eqref{stracchino}
we get \eqref{giuncata3}.
\eproof

\medskip
Take $\epsilon>\e$.
Recalling \eqref{enorme},
 \eqref{trescul0} and \eqref{CEV}
for every $N$ large enough, we get
\begin{equation}\label{trescul2}
|u_0|_{H^{s+2}_N}^2=
\sum_{j\leq N}j^{2}u_{0j}^2
+\sum_{j> N}
j^{2s+4}u_{0j}^2\leq \epsilon^2
\,,\quad 
|v_0|_{\tilde H^{s}_N}^2
=
\sum_{j\leq N}\frac{j^{2}}{\o_j^2}v_{0j}^2
+\sum_{j> N}
\frac{j^{2s+4}}{\o_j^2}v_{0j}^2\leq \epsilon^2\,.
\end{equation}
By \eqref{nduja} and
\eqref{tommaso}
we have that the initial datum $z_0$
satisfies
\begin{equation}\label{frittata}
|z_{0j}|^2=\frac{q_{0j}^2+p_{0j}^2}{2}
=
\frac{\o_j u_{0j}^2+\o_j^{-1}v_{0j}^2}{2}\,.
\end{equation}
Defining the weight $\tw=\tw_N^s$ through
$$
\tw_{N,j}^s:=
\left\{
\begin{array}{ll}
\displaystyle j\o_j^{-1/2}
\ \ \ \ &\mathrm{if}\ \ \ j\leq N\\
\\
\displaystyle j^{s+2}\o_j^{-1/2}\ \ \ \ &\mathrm{if}\ \ \ j>N\,.
\end{array}
\right.
$$
by \eqref{trescul2} and \eqref{frittata}
 we have that
\begin{equation}\label{vulcano}
|z_0|_{\tw^s_N}\leq \epsilon\,.
\end{equation}
Recalling \eqref{peso} and \eqref{mortadella}
we have, for every $N$,
\begin{equation}\label{brunilde}
  \g_r=r\,.
\end{equation}
Setting 
\begin{equation}\label{brunilde2}
r:=\left(\frac56-\frac{1}{80e}\right)^{-1}\epsilon\,,
\end{equation}
by 
\eqref{trescul0}
the hypotheses of Corollary \ref{moderno}
(recall also \eqref{daje})
are satisfied
since
$$
0.08<
\left(\frac56-\frac{1}{80e}\right)
 \frac{5}{11}\sqrt{\frac{\pi}{22e}}\,.
$$
Then
by \eqref{giuncata3}
we get, for every $|t|\leq T_0$,
 \begin{equation}\label{giuncata4}
|z(t)|_\tw\leq 
\left(1+\frac{1}{80e}\right)
\left(\frac56-\frac{1}{80e}\right)^{-1}\epsilon
\end{equation}
Recalling \eqref{peppa}
we get 
$$
|u|_{H^{s+2}_N}\leq \sqrt 2 |z|_{\tw^s_N}
$$
and the first estimate in \eqref{burrata} follows taking $\epsilon$ close enough to $\e$
and $N$ large enough,
 for every $|t|\leq T_0$, by 
\eqref{giuncata4}
since
$$
\left(1+\frac{1}{80e}\right)
\left(\frac56-\frac{1}{80e}\right)^{-1}\sqrt 2
<1.9\,.
$$
The second estimate in \eqref{burrata} follows by the first one 
noting that $1.9\sqrt{\pi/3}<2$.
Let us finally evaluate $T_0$ in \eqref{stoccafisso}.
By 
\eqref{duspaghibis},
 \eqref{brunilde} and \eqref{brunilde2}
we have 
\begin{equation}\label{aerosol2}
\left(\frac{11\g_r }{10}\right)^2
=c_*\epsilon^2\,,\qquad
\eta_r=c_\dag c_*\epsilon^2/\mu\,.
\end{equation}
Then
 by
 \eqref{giuncata2}
 and noting that $\pi e/9>0.948$,
 the expression for
 $T_0$ in \eqref{stoccafisso} follows
 taking $\epsilon$ close enough to $\e$
and $N$ large enough.
\eproof

\bigskip

\noindent
\subsubsection*{Proof of Theorem \ref{patataTER}}

In view of \eqref{stoccafisso},
\eqref{golfetta}, \eqref{bagno}, \eqref{brunilde}
and \eqref{aerosol2}, by \eqref{duspaghibis}
we have
\begin{equation}\label{frattaje}
	 C_r\leq C_r^*:=
	\frac{132\cccp^2}{\pi^2 \mu}
 \left(\frac{11 }{10}\right)^4 r^4
 =\frac{132\cccp^2c_*^2}{\pi^2 \mu}\epsilon^4
 \leq 41.6 \frac{\cccp^2}{ \mu}\epsilon^4
 \,. 
\end{equation}
Recalling \eqref{daje6} and that $\tr=10r/11$
we get
\begin{equation}\label{daje69}
 \tilde\eta_\tr\leq 
\frac{352 e C_{11\tr/10}^*}{7|m|}
=
\frac{352 e C_r^*}{7|m|}
\stackrel{\eqref{frattaje}}
=
\frac{46464 e \cccp^2c_*^2}{7|m|\pi^2 \mu}\epsilon^4
\leq
\frac{ 3221\cccp^2}{|m| \mu}\epsilon^4
=:\tilde\eta_\tr^*
\,.
\end{equation}
Then
condition \eqref{daje6}
is satisfied assuming 
that $\tilde\eta_\tr^*\leq 1/2$, namely
\begin{equation}\label{babele}
\e\leq\sqrt[4]{\frac{|m| \mu}{ 3221\cccp^2}}
\end{equation}
which follows by \eqref{e1}
 taking $\epsilon$ close enough to $\e$
and $N$ large enough.
Recalling \eqref{duspaghibis6}
and 
\eqref{brunilde2}
\begin{eqnarray*}
	 \tilde C_\tr
	 &\leq&
	\frac{1}{44e^2}\tilde\eta_\tr^2+	
	\frac{2}{e}\tilde\eta_\tr
\left(\frac{\cccp}\pi r^2
+C_r^*
\right)+
\frac{16\cccp}{3e^2\pi^3}
\frac{r^6}{\mu^2} 
(1+1 /e)
\nonumber
\\
&\leq&
(31911+98587|m|)
\frac{\cccp^4}{|m|^2 \mu^2}\epsilon^8
+1099\frac{\cccp^3}{|m|\mu}\epsilon^6
\label{cinghiale}
\\
&\stackrel{\eqref{babele}}\leq&
\Big(1+
29.1\frac{\cccp}{|m|\mu}\epsilon^2+
89.8\frac{\cccp}{\mu}\epsilon^2
\Big)
1099\frac{\cccp^3}{|m|\mu}\epsilon^6
\leq \frac{1}{6T_1}\e^6\,,
\label{cinghiale2}
\end{eqnarray*}
 taking $\epsilon$ close enough to $\e$
and $N$ large enough.
We finally apply Lemma \ref{cobra} to 
the Hamiltonian $\mathtt H$
in \eqref{pesce6} with ${\mathcal N}\to
\Lambda + \bar G+\bar R$, $R\to\mathtt R$, $r\to\tr$. Then by \eqref{duspaghibis6} we 
get 
that every solutions of $\dot \z = X_{\mathtt H}$
satisfying
  \[
  \z(0)=\z_0\,, \qquad
  |\z_0|_\tw<\frac56 \tr\,,
 \]
 verify
 \begin{equation}\label{giuncata6}
\Big| |\z(t)|_\tw-|\z_0|_\tw\Big|< \frac{\tr}{6}\,, 
\quad
|\z(t)|_\tw\leq \tr
\qquad
\forall\,  |t|\leq 
\frac{1}{6 \tilde C_\tr}
\,.
\end{equation}

\section{Appendix}\label{app}
		
		\noindent
\textbf{Proof of Proposition
\ref{ham flow}}
	For brevity we set 
	$
	|\cdot|_{r'}:=|\cdot|_{r',\tw}\,,
	$
	for every $r'>0$.
We need the following elementary result 
(see Lemma B.3 of \cite{BMP})
\begin{lem}
 Let $0<r_1<r.$ Let $E$ be a Banach space endowed with the norm $|\cdot|_E$.
 Let $X:B_r \to E$ a vector field satisfying
 $$\sup_{B_r}|X|_E\leq \delta_0\,.$$
 Then the flow $\Phi(u,t)$ of the vector field\footnote{Namely the solution 
 of the equation $\partial_t \Phi(u,t)=X(\Phi(u,t))$ with initial datum
 $\Phi(u,0)=u.$} is well defined for every 
 $$|t|\leq T:=\frac{r-r_1}{\delta_0}$$
 and $u\in B_{r_1}$
 with estimate
 $$
 |\Phi(u,t)-u|_E\leq \delta_0 |t|\,,\qquad
 \forall\, |t|\leq T
 \,.
 $$
\end{lem} 	
	Since
	\[
	\sup_{u\in  B_{r+\rho}(\th_{\mathtt w})}
	|X_S|_{\th_{\mathtt w}}
	\le (r+\rho) |S|_{r+\rho}=\frac{\rho\eta}{4 e}\,,
	\]
	by the above lemma with $E\to \th_{\mathtt w}$, $X\to X_S$,
	$\delta_0\to {\rho\eta}/{4 e},$ $r\to r+\rho,$ $r_1\to r,$	we have that the time $1$-Hamiltonian flow 
	$\Phi^1_S: B_r(\th_{\mathtt w})
	\to B_{r + \rho}(\th_{\mathtt w})$  is well defined
	and analytic, and that
	\[
	\sup_{u\in  B_{r}(\th_{\mathtt w})}
	\abs{\Phi^1_S(u)-u} _{\th_{\mathtt w}}
	\le
	\frac{\rho\eta}{4 e}
	\,.
	\] 
	Estimates \eqref{caio} directly follows by
	\eqref{Lie} and \eqref{brubeck} with  $c_n=1/n!$.
		Let us prove \eqref{brubeck}.
	Fix $n\in\N,$ $n>0$ and set
	$$
	r_i := r +\rho(1 - \frac{i}{n}) \,
	\, ,  \qquad  i = 0,\ldots,n \, .
	$$
	Note that, by the monotonicity of the norm 	\begin{equation}\label{morello}
	|S|_{r_i}\leq |S|_{r+\rho}\,,\qquad
	\forall\,  i = 0,\ldots,n\,.
	\end{equation}
	Noting that
	\begin{equation}\label{dave}
	1+\frac{n r_i}{\rho} \,
	\leq
	n \pa{1+\frac{r}{\rho }}\,,
	\qquad \forall\, i=0,\ldots,n\,,
	\end{equation}
	by using $n$ times \eqref{commXHK}  we have
	\begin{eqnarray*}
		| {H^{(n)}}|_r
		&=& 
		|   \{S, {H^{(n-1)}}\} |_r
		\leq  
		4 (1+\frac{ n r}{\rho})
		|{H^{(n-1)}}|_{r_{n-1}}|
		S|_{r_{n-1}}
		\\
		&\stackrel{\eqref{morello}}\leq&
		|H|_{r+\rho}
		|S|_{r+\rho}^n
		4^n
		\prod_{i=1}^n
		(1+\frac{ n r_i}{\rho})
		\stackrel{\eqref{dave}}\leq
		|H|_{r+\rho}
		\left(
		4n \pa{1+\frac{r}{\rho }}|S|_{r+\rho}
		\right)^n
		\,.
	\end{eqnarray*} 
	Then
	we  get
	\begin{eqnarray*}
		\left|\sum_{n\geq h} c_n {H^{(n)}}\right|_{r} &\leq&
		\sum_{n\geq h} |c_n| |{H^{(n)}}|_{r}
		\leq
		|H|_{r+\rho} \sum_{n\geq h}|c_n| \left(
		4n \pa{1+\frac{r}{\rho }}|S|_{r+\rho}
		\right)^n
		\\
		&\leq &
		\sup_{n\geq h}  \big(|c_n|n^n e^{-n}\big) | H|_{r+\rho}
		 \sum_{n\geq h} \eta^n
		\,,
	\end{eqnarray*}	
	proving \eqref{brubeck}.
\eproof  			

\medskip

\noindent
\textbf{Proof of Lemma \ref{cobra}}
Let us look at the time evolution of $|z(t)|_\tw^2$.  
 By construction and Cauchy-Schwarz 
 inequality
 \begin{eqnarray*}
 2|z(t)|_\tw  \left|\frac{d}{dt} |z(t)|_\tw \right|
 &=&
 \left|\frac{d}{dt} |z(t)|_\tw^2 \right|
 = 
2| {\rm Re} (z,\dot z)_{\th_\tw}|
=
2| {\rm Re} (z,X_R)_{\th_\tw}|
 \\
 &\le&
  2  |z(t)|_\tw |X_{\und R}|_\tw 
\leq
  2 r |z(t)|_\tw |R|_{r,\tw}\
\end{eqnarray*}
 as long as $|z(t)|_\tw\leq r$;
 namely
 \begin{equation}\label{bufala}
\left|\frac{d}{dt} |z(t)|_\tw\right|
 \leq
 r |R|_{r,\tw}
\end{equation}
  as long as $|z(t)|_\tw\leq r.$
Assume by contradiction that there exists
a time\footnote{The case $T_0<0$ is analogous.}  
  $$
  0<T_0<\frac{1}{6|R|_{r,\tw}}
  $$ 
  such that 
  \begin{equation}\label{giuncata5}
\Big| |z(t)|_\tw-|z_0|_\tw\Big|
< \frac{r}{6}\,, \quad
\forall\, 0\leq t<T_0\,,\qquad {\rm but}\ \ 
\Big| |z(T_0)|_\tw-|z_0|_\tw\Big|= \frac{r}{6}
\,.
\end{equation}
  Then
  $$
  |z(t)|_\tw\leq |z_0|_\tw+ \frac{r}{6}
 < 
 r\,
  \qquad
  \quad
\forall\, 0\leq t\leq T_0
  \,.
  $$
  By \eqref{bufala}  we get
  $$
  \Big| |z(T_0)|_\tw-|z_0|_\tw\Big|
  \leq
  r |R|_{r,\tw} T_0
  <\frac{r}{6}\,,
  $$
  which contradicts \eqref{giuncata5},
  proving \eqref{giuncata}. 
  \eproof

\bigskip
\bigskip

\end{document}